\pdfoutput=1
\RequirePackage{ifpdf}
\ifpdf 
\documentclass[pdftex]{sigma}
\else
\documentclass{sigma}
\fi

\numberwithin{equation}{section}

\newtheorem{Theorem}{Theorem}[section]
\newtheorem*{Theorem*}{Theorem}

 { \theoremstyle{definition}

\newtheorem{Remark}[Theorem]{Remark}
\newtheorem{Remarks}[Theorem]{Remarks}
}

\begin{document}
\allowdisplaybreaks

\renewcommand{\thefootnote}{}

\newcommand{\arXivNumber}{2307.03616}

\renewcommand{\PaperNumber}{088}

\FirstPageHeading

\ShortArticleName{A Poincar\'e Formula for Differential Forms and Applications}

\ArticleName{A Poincar\'e Formula for Differential Forms\\ and Applications\footnote{This paper is a~contribution to the Special Issue on Global Analysis on Manifolds in honor of Christian B\"ar for his 60th birthday. The~full collection is available at \href{https://www.emis.de/journals/SIGMA/Baer.html}{https://www.emis.de/journals/SIGMA/Baer.html}}}

\Author{Nicolas GINOUX~$^{\rm a}$, Georges HABIB~$^{\rm ab}$ and Simon RAULOT~$^{\rm c}$}

\AuthorNameForHeading{N.~Ginoux, G.~Habib and S.~Raulot}

\Address{$^{\rm a)}$~Universit\'e de Lorraine, CNRS, IECL, F-57000 Metz, France}
\EmailD{\href{nicolas.ginoux@univ-lorraine.fr}{nicolas.ginoux@univ-lorraine.fr}}
\URLaddressD{\url{https://nicolas-ginoux.perso.math.cnrs.fr/}}

\Address{$^{\rm b)}$~Lebanese University, Faculty of Sciences II, Department of Mathematics, \\
\hphantom{$^{\rm b)}$}~P.O. Box 90656 Fanar-Matn, Lebanon}
\EmailD{\href{ghabib@ul.edu.lb}{ghabib@ul.edu.lb}}
\URLaddressD{\url{https://iecl.univ-lorraine.fr/membre-iecl/habib-georges/}}

\Address{$^{\rm c)}$~Universit\'e de Rouen Normandie, CNRS, Normandie Univ, LMRS UMR 6085,\\
\hphantom{$^{\rm c)}$}~F-76000 Rouen, France}
\EmailD{\href{simon.raulot@univ-rouen.fr}{simon.raulot@univ-rouen.fr}}
\URLaddressD{\url{https://lmrs.univ-rouen.fr/persopage/simon-raulot}}

\ArticleDates{Received July 19, 2023, in final form October 26, 2023; Published online November 08, 2023}

\Abstract{We prove a new general Poincar\'e-type i\-ne\-qua\-li\-ty for dif\-fe\-ren\-tial forms on compact Riemannian manifolds with nonempty boundary. When the boundary is isometrically immersed in Euclidean space, we derive a new i\-ne\-qua\-li\-ty involving mean and scalar curvatures of the boundary only and characterize its limiting case in codimension one. A new Ros-type i\-ne\-qua\-li\-ty for differential forms is also derived assuming the existence of a nonzero parallel form on the manifold.}


\Classification{53C21; 53C24; 58J32; 58J50}

\begin{flushright}
{\it Dedicated to Christian B\"ar for his sixtieth birthday}
\end{flushright}

\renewcommand{\thefootnote}{\arabic{footnote}}
\setcounter{footnote}{0}

\section{Introduction}

This article can be seen as a step in understanding the general effect of the curvature operator of the interior of a ma\-ni\-fold with boundary on the to\-po\-lo\-gy and
the geo\-me\-try of its boundary.
To motivate our precise setting, let us recall some previous
works in this area.
In \cite{ShiTam02}, inspired by questions arising in general relativity such as the
po\-si\-ti\-vi\-ty of the Brown--York quasi-local mass \cite{BY}, Shi and Tam shed light
on how the scalar curvature affects the total mean curvature of the boundary.
Namely, if $\big(M^{n+1},g\big)$ is a compact Riemannian spin ma\-ni\-fold~-- for basics on spinors and Dirac operators, see, e.g., \cite{BHMM,Friedrichbook2000,GinsurveyDirac,LawsonMichelsohn89}~--
with nonnegative scalar curvature such that its boundary can be isometrically
embedded into the Euclidean space as a strictly convex hypersurface, then the
total mean curvature of the boundary cannot be greater than the one of the
Euclidean embedding.
More precisely, if $H$ (resp.~$H_0$) denotes the mean curvature of the boundary $\partial M$ in $M$
(resp. in the Euclidean space) then
\begin{equation}\label{ShiTamintro}
\int_{\partial M}H\,{\rm d}\mu_g\leq \int_{\partial M} H_0 \,{\rm d}\mu_g,
\end{equation}
and the equality is attained if and only if the ma\-ni\-fold $M$ is isometric to a
domain in the Euclidean space.
This result  has a lot of deep and important consequences both in mathematics
and physics.
In particular, although its proof relies on the positive mass theorem (PMT), it
is shown to be actually equivalent to this famous result of ma\-the\-ma\-ti\-cal general
relativity.

In the same spirit and with an alternative method, Hijazi and Montiel
\cite{HijaziMontiel14} showed that an~inequa\-li\-ty similar to (\ref{ShiTamintro}) can be deduced from
a general integral inequa\-li\-ty which holds for any spinor field defined on $\partial M$.
Such i\-ne\-qua\-li\-ties will be referred to as {\it Poincar\'e type i\-ne\-qua\-li\-ties} in
the following.
More precisely, they proved that if $\big(M^{n+1},g\big)$ is a compact Riemannian spin
ma\-ni\-fold with mean convex smooth boundary $\partial M$ (endowed with the induced
Riemannian and spin structures) and if $D$ denotes the Dirac operator acting on
the spinor bundle $\Sigma\partial M$ over $\partial M$, then \looseness=-1
\begin{gather}\label{IntIneq}
\frac{n^2}{4}\int_{\partial M} H|\varphi|^2\,{\rm d}\mu_g\leq\int_{\partial M}\frac{|D\varphi|^2}{H}\,{\rm d}\mu_g
\end{gather}
for all $\varphi\in\Gamma(\Sigma\partial M)$.
When the boundary $\partial M$ can be isometrically immersed into another Riemannian spin
ma\-ni\-fold carrying a parallel spinor $\varphi$ with mean curvature $H_0$, the
restriction of such a spinor field to $\partial M$ provides a solution to the Dirac equation
\[
D\varphi=\frac{n}{2}H_0\varphi\qquad\textrm{and}\qquad |\varphi|=1.
\]
Using this spinor field in the i\-ne\-qua\-li\-ty (\ref{IntIneq}) yields
\begin{equation}\label{eq:1}
\int_{\partial M}H\,{\rm d}\mu_g\leq \int_{\partial M} \frac{H_0^2}{H}\,{\rm d}\mu_g
\end{equation}
and equality is characterized by both immersions having the same
shape o\-pe\-ra\-tors.
Even if this i\-ne\-qua\-li\-ty is a straightforward consequence of (\ref{ShiTamintro}) and
the Cauchy--Schwarz i\-ne\-qua\-li\-ty, it has interesting counterparts.
First, it holds under less stringent conditions since it only requires
mean convexity while strict convexity is needed in the Shi--Tam i\-ne\-qua\-li\-ty.
Moreover, the assumption of existence of an isometric immersion in the Euclidean space is relaxed to allow more general ambient spaces, which include, among others, Calabi--Yau and hyperk\"ahler ma\-ni\-folds.
A further interesting property is that, although the proof of the i\-ne\-qua\-li\-ty
(\ref{eq:1}) does not use the PMT, it still implies this result, at least in the
case $n=2$.

Motivated by those results, Miao and Wang \cite{MiaoWang14} considered the same geo\-me\-tric
set-up as before but without assuming the spin condition.
They showed that i\-ne\-qua\-li\-ties similar to (\ref{eq:1}) could be proved by requiring a
lower bound, say $K$, on the Ricci curvature of the Riemannian ma\-ni\-fold
$\bigl(M^{n+1},g\bigr)$ instead of the scalar curvature.
The condition on the Ricci curvature comes naturally in this context and, in a first step, a Poincar\'e type i\-ne\-qua\-li\-ty \cite[Theorem 2.1]{MiaoWang14} can be established via the Reilly formula and which reads as
\begin{equation}\label{Poincare}
\int_{\partial M}\langle S\nabla f,\nabla f\rangle\,{\rm d}\mu_g\leq \int_{\partial M}\frac{1}{H}(\Delta f-tf)^2
\,{\rm d}\mu_g
\end{equation}
for any smooth function $f$ on $\partial M$,  any constant $t\leq \frac{1}{2}K$, and where $\Delta$ is the Laplace operator acting on functions on $M$ and $S$ is the second fundamental form of $\partial M$ in $M$.
Assuming furthermore the existence of an isometric immersion $X\colon
\partial M\rightarrow\mathbb{R}^m$ of $\partial M$ into some Euclidean space $\mathbb{R}^m$ with
$m\geq n+1$ and using the components $x_j$ of this immersion in (\ref{Poincare})
for all $j=1,\dots,m$, Miao and Wang deduced that
\[
\int_{\partial M}H\,{\rm d}\mu_g\leq \int_{\partial M}\frac{\bigl|\vec{H}_0\bigr|^2}{H}\,{\rm d}\mu_g,
\]
where $\vec{H}_0$ is the mean curvature vector of the immersion $X$.
Note that the e\-xis\-ten\-ce of such an immersion is guaranteed by the Nash
embedding theorem.
As a corollary of that i\-ne\-qua\-li\-ty, Miao and Wang obtain rigidity results for ma\-ni\-folds with boundary and Ricci curvature bounded from
below.
It is also important to note that the i\-ne\-qua\-li\-ty (\ref{Poincare}) has a natural
physical interpretation since, as noticed in \cite{MT,MTX}, it appears in the
second variation of the Wang--Yau quasi-local mass \cite{WY}.

In the present article, we generalize (\ref{Poincare}) to differential forms of arbitrary degree, assuming  suitable but very general curvature conditions on the interior as well as on the boundary of the manifold, see i\-ne\-qua\-li\-ty (\ref{eq:Poincareineqpformsgeneral}) in Theorem \ref{t:Poincareineqpformsgeneral} as well as Theorem \ref{thm:ineomegabiggerk}.
The special case, where~(\ref{eq:Poincareineqpformsgeneral}) is an equality turns out to be very rigid, imposing restrictions on the differential form and the geometry of the underlying manifold.
In view of the Shi--Tam i\-ne\-qua\-li\-ty, we look at the particular case where the boundary is isometrically immersed in Euclidean space and are able to both simplify the i\-ne\-qua\-li\-ty and deduce a rigidity result extending Miao and Wang's one in case the boundary is one-codimensional in some affine subspace, see Theorem \ref{t:PoincareineqpformsboundaryinEuclideanspace}.
When the boundary can be immersed in a round sphere, the corresponding i\-ne\-qua\-li\-ty turns out to be strict, see Theo\-rem \ref{t:MwithdelMinSn+m}.
In a next step, we adapt to the differential-form-framework the celebrated Ros i\-ne\-qua\-li\-ty \cite[Theorem 1]{Ros87} involving the integral of the inverse of the mean curvature over the boundary.
Assuming the existence of a nonzero parallel form on the manifold, a new i\-ne\-qua\-li\-ty relating the integral of the inverse of some $\sigma_p$-curvature on the boundary with the volume of the manifold can be deduced from the so-called Reilly formula, see Theorem \ref{t:Rosgeneral}.

The article is structured as follows.
After preliminaries about basic formulae and notations in Section \ref{s:prelim}, the main Poincar\'e-type i\-ne\-qua\-li\-ty (\ref{eq:Poincareineqpformsgeneral}) is presented and proved in Section \ref{s:Poincareineqpformsgeneral}.
When the boundary can be isometrically immersed in Euclidean space, the i\-ne\-qua\-li\-ty can be simplified as we mentioned above, while the case where the interior curvature condition is relaxed is presented in Theorem \ref{thm:ineomegabiggerk}.
In Section \ref{s:Rosdiffforms}, a differential-form-version of the Ros i\-ne\-qua\-li\-ty is established.

\section{Preliminaries and notations}\label{s:prelim}

In this section, we briefly introduce the geometric setting and fix the notations of this paper.

Let $\big(M^{n+1},g\big)$ be a compact oriented $(n+1)$- dimensional Riemannian ma\-ni\-fold with smooth nonempty boundary $\partial M$ and let $\iota\colon\partial M\to M$ be the inclusion map.
Let ${\rm d}\mu_g$ be the Riemannian measure induced by $g$ on both $M$ and $\partial M$.
In the following, any metric -- being $g$ on $TM$ or any metric induced by $g$ on further bundles -- will be denoted by $\langle\cdot\,,\cdot\rangle$, with associated pointwise norm $|\cdot|$.
Let $\nu$ denote the inward unit normal along $\partial M$ and $S:=-\nabla\nu$ be the associated Weingarten endomorphism-field, where $\nabla$ denotes the Levi-Civita connection on $TM$.
Let $H:=\frac{1}{n}\operatorname{tr}(S)$ denote the mean curvature of $\partial M$ in $M$.
For any integer $p\in\{0,\dots,n+1\}$, let $\Omega^p(M):=\Gamma(\Lambda^pT^*M)$ be the space of differential forms of degree $p$ on $M$, that is the space of sections of the exterior bundle $\Lambda^p T^*M\to M$.
Let $\star\colon\Lambda^pT^*M\to\Lambda^{n+1-p}T^*M$ denote the pointwise Hodge star  operator.
For any $(1,1)$-tensor $A$ and $p$-form $\omega$, let $A^{[p]}\omega$ be the $p$-form that is pointwise defined by the following identity: for all tangent vectors $X_1,\dots,X_p$,
\[\bigl(A^{[p]}\omega\bigr)(X_1,\dots,X_p):=\sum_{j=1}^p\omega(X_1,\dots,AX_j,\dots,X_p).\]
In the particular case, where $A=S$, we denote for each $k\in\{1,\dots,n\}$ by $\sigma_k$ the pointwise $k$-curvature of $\partial M$, that is the sum of the $k$ smallest principal curvatures (i.e., the eigenvalues of~$S$) of $\partial M$.
Note that, since the eigenvalues of $S$ are the sums of exactly $p$ among the $n$~principal curvatures,
\[\bigl\langle S^{[p]}\omega,\omega\bigr\rangle\geq\sigma_p|\omega|^2\]
for any $\omega\in\Lambda^pT^*\partial M$, with equality if and only if $S^{[p]}\omega=\sigma_p\omega$.
Moreover, for all $1\leq p\leq q\leq n$, we have $\frac{\sigma_p}{p}\leq\frac{\sigma_q}{q}$, with equality when $p<q$ if and only if the $q$ smallest principal curvatures are equal.

Let ${\rm d}$ (resp.\ $\delta$) denote the exterior derivative (resp.\ codifferential) on $p$-forms and $\nabla$ be the covariant derivative induced by $\nabla$ on $\Lambda^pT^*M$.
Recall the so-called \emph{Reilly formula} \cite[Theorem~3]{RaulotSavo11} for differential $p$-forms with $p\geq 1$: for any $\omega\in\Omega^p(M)$
\begin{gather}
\int_M\bigl(|{\rm d}\omega|^2+|\delta\omega|^2-|\nabla\omega|^2-\bigl\langle W^{[p]}\omega,\omega\bigr\rangle\bigr)\,{\rm d}\mu_g \nonumber\\
\qquad{} =\int_{\partial M}\bigl(2\bigl\langle\nu\lrcorner\,\omega,\delta^{\partial M}(\iota^*\omega)\bigr\rangle+\bigl\langle S^{[p]}\iota^*\omega,\iota^*\omega\bigr\rangle+\bigl\langle S^{[n+1-p]}\iota^*(\star\omega),\iota^*(\star\omega)\bigr\rangle\bigr)\,{\rm d}\mu_g,\label{eq:Reillyformula}
\end{gather}
where $\delta^{\partial M}$ denotes the codifferential on $\partial M$ and $W^{[p]}$ the curvature term involved in the Weitzenb\"ock formula for $p$ forms: denoting by $\Delta:={\rm d}\delta+\delta {\rm d}$ the Hodge Laplace operator on $p$-forms,
\[\Delta\omega=\nabla^*\nabla\omega+W^{[p]}\omega\]
for any $p$-form $\omega$.  By convention, we let $W^{[p]}=0$ and $S^{[p]}=0$ for all $p\leq 0$.
Note that \cite[Theorem 3]{RaulotSavo11}
\begin{equation*}
\bigl\langle S^{[n+1-p]}\iota^*(\star\omega),\iota^*(\star\omega)\bigr\rangle=nH|\nu\lrcorner\,\omega|^2-\bigl\langle S^{[p-1]}(\nu\lrcorner\,\omega),\nu\lrcorner\,\omega\bigr\rangle
\end{equation*}
for any $p$-form $\omega$.

\section[A Poincar\'e-type inequality for p-forms]{A Poincar\'e-type inequality for $\boldsymbol p$-forms}\label{s:Poincareineqpformsgeneral}

We first prove a generalized version of the integral i\-ne\-qua\-li\-ty for functions obtained by Miao and Wang in \cite[equation~(1.3)]{MiaoWang14}.

\begin{Theorem}\label{t:Poincareineqpformsgeneral}
Let $\bigl(M^{n+1},g\bigr)$ be any compact oriented  $(n+1)$-dimensional Riemannian manifold with smooth nonempty boundary $\partial M$.
Fix $p\in\{1,\dots,n\}$ and assume that $W^{[p]}\geq0$ on $M$ as well as $\sigma_{n+1-p}>0$ along $\partial M$.
Then, for any exact $p$-form $\omega$ on $\partial M$,  we have
\begin{equation}\label{eq:Poincareineqpformsgeneral}
\int_{\partial M}\frac{\bigl|\delta^{\partial M}\omega\bigr|^2}{\sigma_{n+1-p}}\,{\rm d}\mu_g\geq\int_{\partial M}\langle S^{[p]}\omega,\omega\rangle\,{\rm d}\mu_g.
\end{equation}
Moreover, if {\rm (\ref{eq:Poincareineqpformsgeneral})} is an equality for some non-zero $\omega$, then for the $p$-form $\hat{\omega}$ on $M$ satisfying ${\rm d}\hat{\omega}=0=\delta\hat{\omega}$ on $M$ as well as $\iota^*\hat{\omega}=\omega$ on $\partial M$, the identities $\delta^{\partial M}\omega=-\sigma_{n+1-p}\nu\lrcorner\,\hat{\omega}$, $S^{[p-1]}(\nu\lrcorner\,\hat{\omega})=(nH-\sigma_{n+1-p})\nu\lrcorner\,\hat{\omega}$ hold along $\partial M$ and the $p$-form $\hat{\omega}$ must be parallel -- hence $W^{[p]}\hat{\omega}=0$ must hold -- on $M$.
\end{Theorem}

\begin{proof}
The proof mainly follows that of \cite[Theorem 5]{RaulotSavo11}.
Let $\omega$ be an exact $p$-form on $\partial M$, that is $\omega={\rm d}^{\partial M}\alpha$ for some $\alpha\in\Omega^{p-1}(\partial M)$.
By \cite[Theorem 2]{DuffSpencer52}, there exists a $(p-1)$-form $\hat{\alpha}$ on~$M$ such that $\delta {\rm d}\hat{\alpha}=0$ on $M$ with $\iota^*\hat{\alpha}=\alpha$ along $\partial M$. Let $\hat{\omega}:={\rm d}\hat{\alpha}\in\Omega^p(M)$, then $\hat{\omega}$ satisfies ${\rm d}\hat{\omega}=0$, $\delta\hat{\omega}=\delta {\rm d}\hat{\alpha}=0$ on~$M$ with $\iota^*\hat{\omega}=\iota^*({\rm d}\hat{\alpha})=\omega$ along $\partial M$. Actually such a~$p$-form~$\hat{\omega}$ on~$M$ with ${\rm d}\hat{\omega}=0=\delta\hat{\omega}$ on $M$ as well as $\iota^*\hat{\omega}=\omega$ along $\partial M$ is uniquely determined by $\omega$. Namely, because of $W^{[p]}\geq0$ by assumption and the property $\star W^{[p]}\star^{-1}=W^{[n+1-p]}$, we know that $W^{[n+1-p]}\geq0$. Together with the assumption $\sigma_{n+1-p}>0$ we can deduce from \cite[Theorem~4]{RaulotSavo11} that $H_{\rm abs}^{n+1-p}(M)=0$  holds, where
\[H_{\rm abs}^{k}(M) :=\bigl\{\alpha\in\Omega^k(M) \mid {\rm d}\alpha=0=\delta\alpha\textrm{ and }\nu\lrcorner\,\alpha=0\bigr\}\]
is, for every $0\leq k\leq n+1$, the $k^{\textrm{th}}$ absolute de Rham cohomology group. By Poincar\'e duality, $H_{\rm abs}^{k}(M)\cong H_{\rm rel}^{n+1-k}(M)$,  where
\[H_{\rm rel}^{k}(M):=\bigl\{\alpha\in\Omega^k(M) \mid {\rm d}\alpha=0=\delta\alpha\textrm{ and }\iota^*\alpha=0\bigr\}\]
is the $k^{\textrm{th}}$ relative de Rham cohomology group. Therefore, $H_{\rm rel}^{p}(M)=0$, so that $\hat{\omega}$ is uniquely determined by $\omega$.

We apply identity (\ref{eq:Reillyformula}) to $\hat{\omega}$: since ${\rm d}\hat{\omega}=\delta\hat{\omega}=0$ by construction of $\hat{\omega}$ and $\bigl\langle W^{[p]}\hat{\omega},\hat{\omega}\bigr\rangle\geq0$ by assumption,  as well as $|\nabla\hat{\omega}|^2\geq 0$, we have from (\ref{eq:Reillyformula}) that
\[0\geq\int_{\partial M}\bigl(2\bigl\langle\nu\lrcorner\,\hat{\omega},\delta^{\partial M}\omega\bigr\rangle+\bigl\langle S^{[p]}\omega,\omega\bigr\rangle+\bigl\langle S^{[n+1-p]}\iota^*(\star\hat{\omega}),\iota^*(\star\hat{\omega})\bigr\rangle\bigr)\,{\rm d}\mu_g.\]
By definition of the $k$-curvatures and the identity $\iota^*(\star\hat{\omega})=(-1)^p\nu\lrcorner\star(\nu\lrcorner\,\hat{\omega})=\star_{\partial M}(\nu\lrcorner\,\hat{\omega})$, we have\looseness=-1
\begin{equation}\label{eq:sn1p}
\bigl\langle S^{[n+1-p]}\iota^*(\star\hat{\omega}),\iota^*(\star\hat{\omega})\bigr\rangle\geq\sigma_{n+1-p}|\iota^*(\star\hat{\omega})|^2=\sigma_{n+1-p}|\nu\lrcorner\,\hat{\omega}|^2
\end{equation}
on $\partial M$. Moreover, because $\sigma_{n+1-p}$ is assumed to be positive along $\partial M$,  we write
\begin{gather*}
2\bigl\langle\nu\lrcorner\,\hat{\omega},\delta^{\partial M}\omega\bigr\rangle=\Biggl|\sigma_{n+1-p}^{1/2}\nu\lrcorner\,\hat{\omega}+\frac{1}{\sigma_{n+1-p}^{1/2}}\delta^{\partial M}\omega\Biggr|^2-\sigma_{n+1-p}|\nu\lrcorner\,\hat{\omega}|^2-\frac{|\delta^{\partial M}\omega|^2}{\sigma_{n+1-p}}\\
\hphantom{2\bigl\langle\nu\lrcorner\,\hat{\omega},\delta^{\partial M}\omega\bigr\rangle}
\geq-\sigma_{n+1-p}|\nu\lrcorner\,\hat{\omega}|^2-\frac{|\delta^{\partial M}\omega|^2}{\sigma_{n+1-p}},
\end{gather*}
so that
\[0\geq \int_{\partial M}\left(-\frac{|\delta^{\partial M}\omega|^2}{\sigma_{n+1-p}}+\bigl\langle S^{[p]}\omega,\omega\bigr\rangle\right)\,{\rm d}\mu_g\]
which is inequality   (\ref{eq:Poincareineqpformsgeneral}).

Assume now (\ref{eq:Poincareineqpformsgeneral}) to be an equality for some non-zero $\omega$. Let $\hat{\omega}$ be the $p$-form on $M$ such that ${\rm d}\hat{\omega}=0=\delta\hat{\omega}$ on $M$ and $\iota^*\hat{\omega}=\omega$ on $\partial M$; as mentioned above.
By the sequence of pointwise i\-ne\-qua\-li\-ties used in the above proof of i\-ne\-qua\-li\-ty (\ref{eq:Poincareineqpformsgeneral}), we can deduce that $\nabla\hat{\omega}=0$ on $M$, that is,~$\hat{\omega}$ is parallel on $M$  and, furthermore, $\delta^{\partial M}\omega=-\sigma_{n+1-p}\nu\lrcorner\,\hat{\omega}$ and $S^{[n+1-p]}\iota^*(\star\hat{\omega})=\sigma_{n+1-p}\iota^*(\star\hat{\omega})$ must hold on $\partial M$.  As a straightforward consequence, $W^{[p]}\hat{\omega}=0$ on $M$. Thanks to the identity $S^{[n+1-p]}(\star_{\partial M}\alpha)=-\star_{\partial M}S^{[p-1]}\alpha+nH\star_{\partial M}\alpha$ which is valid pointwise for all $(p-1)$-forms $\alpha$, the latter is equivalent to $S^{[p-1]}(\nu\lrcorner\,\hat{\omega})=(nH-\sigma_{n+1-p})\nu\lrcorner\,\hat{\omega}$ along $\partial M$. This concludes the proof of~Theorem \ref{t:Poincareineqpformsgeneral}.
\end{proof}

\begin{Remark}\label{r:LambdapTMtrivial} Note that, for $1\leq p\leq n$, if the bundle $\Lambda^p T^*M\to M$ is trivialized by such parallel $p$-forms $\hat{\omega}$, then on the one hand the manifold $M$ is flat and, on the other hand, $S^{[p-1]}=(nH-\sigma_{n+1-p})\cdot\mathrm{Id}$ must hold pointwise on $\Lambda^{p-1}T^*\partial M$.
The first statement is a consequence from the fact that the curvature of the manifold $M$ vanishes as soon as that of $\Lambda^p T^*M$ does.
The second statement comes from the map $\Lambda^p T_x^*M\to \Lambda^p T^*_x\partial M\oplus \Lambda^{p-1} T^*_x\partial M$, $\omega\mapsto (\iota^*\omega,\nu\wedge(\nu\lrcorner\,\omega))$, being an isomorphism at any $x\in \partial M$.  In case $p\geq2$ (for $p=1$ that identity is trivial because of $S^{[0]}=0$ by convention and $\sigma_{n}=nH$ by definition), this shows that $\iota\colon\partial M\to M$ must be \emph{totally umbilical}.
\end{Remark}

Next we turn to the case where the boundary of $M$ is assumed to be isometrically immersed in some Euclidean space.

\begin{Theorem}\label{t:PoincareineqpformsboundaryinEuclideanspace}
Let $\bigl(M^{n+1},g\bigr)$ be any compact oriented $(n+1)$-dimensional Riemannian manifold with smooth nonempty boundary $\partial M$.
Assume that $W^{[p]}\geq0$ on $M$ as well as $\sigma_{n+1-p}>0$ along $\partial M$ for a given $p\in\{1,\dots,n\}$.
Assume also that there exists an isometric immersion $\iota_0\colon\partial M\to\mathbb{R}^{n+m}$ with mean curvature vector $H_0$. Then
\begin{equation}\label{eq:PoincareineqpformsboundaryinEuclideanspace}
\int_{\partial M}\frac{n|H_0|^2-\frac{p-1}{n(n-1)}{\rm Scal}^{\partial M}}{\sigma_{n+1-p}}\,{\rm d}\mu_g\geq\int_{\partial M}H\,{\rm d}\mu_g,
\end{equation}
where  ${\rm Scal}^{\partial M}$  denotes the scalar curvature of $\partial M$. Equality holds in {\rm (\ref{eq:PoincareineqpformsboundaryinEuclideanspace})} when $M$ is the $(n+1)$-dimensional flat disk $\mathbb{D}^{n+1}$ standardly embedded in some $(n+1)$-dimensional affine subspace of~$\mathbb{R}^{n+m}$. Conversely, if {\rm (\ref{eq:PoincareineqpformsboundaryinEuclideanspace})} is an equality, $\partial M$ is connected, $p\geq 2$ and $\iota_0(\partial M)$ is contained in some $(n+1)$-dimensional affine subspace of $\mathbb{R}^{n+m}$, then, up to rescaling the metrics on $M$ and on $\mathbb{R}^{n+m}$, the manifold $M$ is isometric to the $(n+1)$-dimensional flat disk $\mathbb{D}^{n+1}$ standardly embedded in that subspace.
\end{Theorem}

\begin{proof}
We take the standard coordinates $(x_1,\dots,x_{n+m})$ on $\mathbb{R}^{n+m}$ and, for any $i_1,\dots,i_p\in\{1,\dots,n+m\}$, we denote by ${\rm d}x_I:={\rm d}x_{i_1}\wedge\dots\wedge {\rm d}x_{i_p}\in\Lambda^p(\mathbb{R}^{n+m})^*$ and by $\omega_I:=\iota_0^* {\rm d}x_I\in\Omega^p(\partial M)$. Note that, since ${\rm d}x_I$ is exact, so is $\omega_I$.
Replacing $\omega$ by $\omega_I$ in (\ref{eq:Poincareineqpformsgeneral}), we obtain
\begin{equation}\label{eq:PoincareomegaI}
\int_{\partial M}\frac{\bigl|\delta^{\partial M}\omega_I\bigr|^2}{\sigma_{n+1-p}}\,{\rm d}\mu_g\geq\int_{\partial M}\bigl\langle S^{[p]}\omega_I,\omega_I\bigr\rangle\,{\rm d}\mu_g.
\end{equation}
We now want to deduce from (\ref{eq:PoincareomegaI}) a more explicit i\-ne\-qua\-li\-ty. For this, we sum (\ref{eq:PoincareomegaI}) over $I$, meaning that we compute the sum $\sum_{i_1,\dots,i_p=1}^{n+m}$ of both sides; mind that the indices $i_1,\dots,i_p$ vary independently and hence are repeated. On the one hand, by \cite[Lemma 2.2]{Savo05},
\begin{equation}\label{eq:sumdeltaomegaI2}
\sum_I\bigl|\delta^{\partial M}\omega_I\bigr|^2=p!\binom{n}{p}p\left(n|H_0|^2-\frac{p-1}{n(n-1)}{\rm Scal}^{\partial M}\right).
\end{equation}
On the other hand, we use the pointwise identity
\[\sum_{j}{\rm e}_j\wedge({\rm e}_j\lrcorner\,\alpha)=k\alpha\]
which is valid for any $k$-form $\alpha$ and any pointwise o.n.b.\ $({\rm e}_j)_j$ of $TM$ or $T\partial M$.
As a straightforward consequence,
\begin{equation*}
\sum_{j=1}^n\langle {\rm e}_j\wedge\alpha,{\rm e}_j\wedge\beta\rangle=(n-k)\langle\alpha,\beta\rangle
\end{equation*}
for any $k$-forms $\alpha$, $\beta$ on an $n$-dimensional space. Thus we may compute the sum of the r.h.s.\ of~(\ref{eq:PoincareomegaI}) as follows:  Let $A\colon T\partial M\to T\partial M$ be any symmetric endomorphism of $\partial M$  and denoting by ${\rm d}x_i^T:=\iota_0^*{\rm d}x_i$ and by $({\rm e}_j)_{1\leq j\leq n}$ a pointwise o.n.b.\ of $T\partial M$, we may write
\begin{align}
\sum_I\bigl\langle A^{[p]}\omega_I,\omega_I\bigr\rangle&= \sum_{i_1,\dots,i_p=1}^{n+m}\bigl\langle A^{[p]}({\rm d}x_{i_1}^T\wedge\dots\wedge {\rm d}x_{i_p}^T),{\rm d}x_{i_1}^T\wedge\dots\wedge {\rm d}x_{i_p}^T\bigr\rangle\nonumber\\
&= \sum_{i_1=1}^{n+m}\sum_{i_2,\dots,i_p=1}^{n+m}\sum_{j,k=1}^n\bigl\langle {\rm d}x_{i_1}^T,{\rm e}_j\bigr\rangle\bigl\langle {\rm d}x_{i_1}^T,{\rm e}_k\bigr\rangle\bigl\langle A^{[p]}({\rm e}_j\wedge {\rm d}x_{i_2,\dots,i_p}^T),{\rm e}_k\wedge {\rm d}x_{i_2,\dots,i_p}^T\bigr\rangle\nonumber\\
&= \sum_{j,k=1}^n\!\underbrace{\left(\sum_{i_1=1}^{n+m}\!\bigl\langle {\rm d}x_{i_1},{\rm e}_j\bigr\rangle\bigl\langle {\rm d}x_{i_1},{\rm e}_k\bigr\rangle\!\right)\!}_{\delta_{jk}}\sum_{i_2,\dots,i_p=1}^{n+m}\!\!\bigl\langle A^{[p]}({\rm e}_j\wedge {\rm d}x_{i_2,\dots,i_p}^T),{\rm e}_k\wedge {\rm d}x_{i_2,\dots,i_p}^T\bigr\rangle\nonumber\\
&= \sum_{j=1}^n\sum_{i_2,\dots,i_p=1}^{n+m}\bigl\langle A^{[p]}({\rm e}_j\wedge {\rm d}x_{i_2,\dots,i_p}^T),{\rm e}_j\wedge {\rm d}x_{i_2,\dots,i_p}^T\bigr\rangle\nonumber\\
&= \sum_{j_1,\dots,j_p=1}^n\bigl\langle A^{[p]}({\rm e}_{j_1}\wedge\dots\wedge {\rm e}_{j_p}),{\rm e}_{j_1}\wedge\dots\wedge {\rm e}_{j_p}\bigr\rangle\nonumber\\
&= p\sum_{j_1,\dots,j_p=1}^n\bigl\langle A{\rm e}_{j_1}\wedge {\rm e}_{j_2}\wedge\dots\wedge {\rm e}_{j_p},{\rm e}_{j_1}\wedge\dots\wedge {\rm e}_{j_p}\bigr\rangle\nonumber\\
&= p(n-p+1)\sum_{j_1,\dots,j_{p-1}=1}^n\bigl\langle A{\rm e}_{j_1}\wedge {\rm e}_{j_2}\wedge\dots\wedge {\rm e}_{j_{p-1}},{\rm e}_{j_1}\wedge\dots\wedge {\rm e}_{j_{p-1}}\bigr\rangle\nonumber\\
&= p(n-p+1)\cdots(n-1)\sum_{j_1=1}^n\bigl\langle A{\rm e}_{j_1},{\rm e}_{j_1}\bigr\rangle\nonumber\\
&= p!\binom{n}{p}\frac{p}{n}{\rm tr}(A).\label{eq:traceap}
\end{align}
When $A=S$ is the associated Weingarten endomorphism-field, we get that
\begin{equation}\label{eq:sumISpomegaI}
\sum_I\bigl\langle S^{[p]}\omega_I,\omega_I\bigr\rangle=p!\binom{n}{p}pH.
\end{equation}
Integrating both (\ref{eq:sumdeltaomegaI2}) (after dividing by $\sigma_{n+1-p}$) and (\ref{eq:sumISpomegaI}) over $\partial M$, we deduce i\-ne\-qua\-li\-ty (\ref{eq:PoincareineqpformsboundaryinEuclideanspace}) from (\ref{eq:PoincareomegaI}).

If $M=\mathbb{D}^{n+1}$ is the $(n+1)$-dimensional flat disk standardly embedded in some $(n+1)$-dimensional affine subspace of $\mathbb{R}^{n+m}$, then $|H_0|=1=H$, $\sigma_{n+1-p}=n+1-p$ and  ${\rm Scal}^{\partial M}=n(n-1)$  along $\partial M=\mathbb{S}^n$ (the $n$-dimensional round sphere of sectional curvature $1$), so that (\ref{eq:PoincareineqpformsboundaryinEuclideanspace})  is an equality.

Conversely, if (\ref{eq:PoincareineqpformsboundaryinEuclideanspace}) is an equality, then for any tuple $I$, (\ref{eq:PoincareomegaI}) must be an equality.
For any $I$, we denote by $\hat{\omega}_I$ the $p$-form on $M$ such that ${\rm d}\hat{\omega}_I=0=\delta\hat{\omega}_I$ on $M$ and $\iota^*\hat{\omega}_I=\omega_I$ along $\partial M$; the existence and uniqueness of $\hat{\omega}_I$ is guaranteed by \cite[Theorem 2]{DuffSpencer52} and the vanishing of  $H_{\rm rel}^p(M)$,  see proof of Theorem \ref{t:Poincareineqpformsgeneral} above. Then Theorem \ref{t:Poincareineqpformsgeneral} implies that, for any $I$, the $p$-form $\hat{\omega}_I$ is parallel on $M$, that $\delta^{\partial M}\omega_I=-\sigma_{n+1-p}\nu\lrcorner\,\hat{\omega}_I$ and that $S^{[p-1]}(\nu\lrcorner\,\hat{\omega}_I)=(nH-\sigma_{n+1-p})\nu\lrcorner\,\hat{\omega}_I$ hold along $\partial M$. By (\ref{eq:sumdeltaomegaI2}), this implies
\begin{gather}\label{eq:sumdeltaomegaI2bis}
\sigma_{n+1-p}\sum_I|\nu\lrcorner\,\hat{\omega}_I|^2=\sum_I\frac{\bigl|\delta^{\partial M}\omega_I\bigr|^2}{\sigma_{n+1-p}}
=p!\binom{n}{p}p\left(\frac{n|H_0|^2-\frac{p-1}{n(n-1)}{\rm Scal}^{\partial M}}{\sigma_{n+1-p}}\right).
\end{gather}
Next we show that $\sum_I|\nu\lrcorner\,\hat{\omega}_I|^2$ is constant along $\partial M$. Differentiating along any $X\in T\partial M$ and using the fact that $\hat{\omega}_I$ is parallel on $M$, we have
\begin{gather*}
X\biggl(\frac{1}{2}\sum_I|\nu\lrcorner\,\hat{\omega}_I|^2\biggr)
=\sum_I\bigl\langle\nabla_X(\nu\lrcorner\,\hat{\omega}_I),\nu\lrcorner\,\hat{\omega}_I\bigr\rangle
=-\sum_I\langle SX\lrcorner\,\hat{\omega}_I,\nu\lrcorner\,\hat{\omega}_I\rangle\\
\qquad{}=-\sum_I\bigl\langle \iota^*(SX\lrcorner\,\hat{\omega}_I),\nu\lrcorner\,\hat{\omega}_I\bigr\rangle
=-\sum_I\langle SX\lrcorner\,\omega_I,\nu\lrcorner\,\hat{\omega}_I\rangle
=-\sum_I\bigl\langle (SX\lrcorner\,{\rm d}x_I)^T,\nu\lrcorner\,\hat{\omega}_I\bigr\rangle.
\end{gather*}

In the following, we will prove that the last sum vanishes.
For this, we will express the $(p-1)$-form $\nu\lrcorner\,\hat{\omega}_I$  in terms of the data of the immersion $\iota_0\colon\partial M\to\mathbb{R}^{n+m}$. That identity will be used at several places in the proof. We denote by $\mathbb{I}$ the second fundamental form of $\iota_0$. Using \cite[equation~(3.3)]{Savo05}, we have that, for each $p$-tuple $I$,
\begin{align}
\delta^{\partial M}\omega_I={} &\sum_{k=1}^p(-1)^{k+1}\mathbb{I}_{{\rm d}x_{i_k}^\perp}^{[p-1]}\bigl({\rm d}x_{i_1}^T\wedge \dots \wedge\widehat{{\rm d}x_{i_k}^T}\wedge\dots \wedge {\rm d}x_{i_p}^T\bigr)\nonumber\\
&{}-n \sum_{k=1}^p(-1)^{k+1}\bigl\langle H_0, {\rm d}x_{i_k}^\perp\bigr\rangle\,{\rm d}x_{i_1}^T\wedge\dots\wedge\widehat{{\rm d}x_{i_k}^T}\wedge\dots \wedge {\rm d}x_{i_p}^T\nonumber\\
 ={} &\sum_{k=1}^p\sum_{a=1}^m(-1)^{k+1}\langle {\rm d}x_{i_k},\nu_a\rangle\mathbb{I}_{\nu_a}^{[p-1]}\bigl({\rm d}x_{i_1}^T\wedge \dots \wedge\widehat{{\rm d}x_{i_k}^T}\wedge\dots \wedge {\rm d}x_{i_p}^T\bigr)\nonumber\\
 &{}-n \sum_{k=1}^p(-1)^{k+1}\langle H_0, {\rm d}x_{i_k}\rangle\,{\rm d}x_{i_1}^T\wedge\dots\wedge\widehat{{\rm d}x_{i_k}^T}\wedge\dots \wedge {\rm d}x_{i_p}^T\nonumber\\
 ={}&\sum_{a=1}^m\mathbb{I}_{\nu_a}^{[p-1]}\bigl((\nu_a\lrcorner\,{\rm d}x_I)^T\bigr)-(nH_0\lrcorner\,{\rm d}x_I)^T,\label{eq:deltaomegaISavo}
\end{align}
where $(\nu_a)_{1\leq a\leq m}$ is a pointwise o.n.b.\ of $T^\perp \partial M$ seen as a subspace of $\mathbb{R}^{n+m}$. Because of $\delta^{\partial M}\omega_I=-\sigma_{n+1-p}\nu\lrcorner\,\hat{\omega}_I$, we obtain
\begin{equation}\label{eq:nuinterioromega}
\nu\lrcorner\,\hat{\omega}_I=-\frac{1}{\sigma_{n+1-p}}\cdot\biggl(\sum_{a=1}^m\mathbb{I}_{\nu_a}^{[p-1]}\bigl((\nu_a\lrcorner\,{\rm d}x_I)^T\bigr)-(nH_0\lrcorner\,{\rm d}x_I)^T\biggr).
\end{equation}
Therefore, in order to show that $\sum_I\bigl\langle (SX\lrcorner\,{\rm d}x_I)^T,\nu\lrcorner\,\hat{\omega}_I\bigr\rangle=0$, it is sufficient to show that both sums $\sum_I\bigl\langle (SX\lrcorner\,{\rm d}x_I)^T,\sum_{a=1}^m\mathbb{I}_{\nu_a}^{[p-1]}\bigl((\nu_a\lrcorner\,{\rm d}x_I)^T\bigr)\bigr\rangle$ and $\sum_I\bigl\langle (SX\lrcorner\,{\rm d}x_I)^T,(nH_0\lrcorner\,{\rm d}x_I)^T\bigr\rangle$ vanish.
Let us make the computation for the first sum, the second can be done in the same way.
We denote by ${\rm e}_J={\rm e}_{j_1}\wedge\dots \wedge {\rm e}_{j_{p-1}}$ with $j_1<j_2<\dots<j_{p-1}$ the orthonormal frame of $\Lambda^{p-1}T_x^*\partial M$ induced by a local orthonormal frame $\{{\rm e}_1,\dots,{\rm e}_n\}$ of $T\partial M$.
For any $a=1,\dots,m$,
\begin{gather*}
\sum_I\bigl\langle (SX\lrcorner\,{\rm d}x_I)^T,\mathbb{I}_{\nu_a}^{[p-1]}\bigl((\nu_a\lrcorner\,{\rm d}x_I)^T\bigr)\bigr\rangle=\sum_{I,J}\bigl\langle (SX\lrcorner\,{\rm d}x_I)^T,{\rm e}_J\bigr\rangle \bigl\langle \mathbb{I}_{\nu_a}^{[p-1]}\bigl((\nu_a\lrcorner\,{\rm d}x_I)^T\bigr),{\rm e}_J\bigr\rangle\\
\hphantom{}=\sum_{I,J}\bigl\langle SX\lrcorner\,{\rm d}x_I,{\rm e}_J\bigr\rangle \bigl\langle \nu_a\lrcorner\,{\rm d}x_I,\mathbb{I}_{\nu_a}^{[p-1]}({\rm e}_J)\bigr\rangle
=\sum_{I,J}\bigl\langle {\rm d}x_I,SX\wedge {\rm e}_J\bigr\rangle \bigl\langle {\rm d}x_I,\nu_a\wedge \mathbb{I}_{\nu_a}^{[p-1]}({\rm e}_J)\bigr\rangle\\
\hphantom{}=\sum_{J}\bigl\langle SX\wedge {\rm e}_J,\nu_a\wedge \mathbb{I}_{\nu_a}^{[p-1]}({\rm e}_J)\bigr\rangle
=0,
\end{gather*}
because of $\nu_a\perp T\partial M$.
Similarly, $\sum_I\bigl\langle (SX\lrcorner\,{\rm d}x_I)^T,(nH_0\lrcorner\,{\rm d}x_I)^T\bigr\rangle=0$ by $H_0\perp T\partial M$.

If $\partial M$ is connected, which will be assumed from now on, then $\sum_I|\nu\lrcorner\,\hat{\omega}_I|^2$ is constant along~$\partial M$. By (\ref{eq:sumdeltaomegaI2bis}) and the assumption that (\ref{eq:PoincareineqpformsboundaryinEuclideanspace}) is an equality, this constant is given by
\begin{align*}
\sum_I|\nu\lrcorner\,\hat{\omega}_I|^2\int_{\partial M}\sigma_{n+1-p}\,{\rm d}\mu_g&= p!\binom{n}{p}p\int_{\partial M}\frac{n|H_0|^2-\frac{p-1}{n(n-1)}{\rm Scal}^{\partial M}}{\sigma_{n+1-p}}\,{\rm d}\mu_g\\
&= p!\binom{n}{p}p\int_{\partial M}H\,{\rm d}\mu_g,
\end{align*}
that is
\begin{equation}\label{eq:constantsumnuintomegaI2}
\sum_I|\nu\lrcorner\,\hat{\omega}_I|^2=p!\binom{n}{p}p\frac{\int_{\partial M}H\,{\rm d}\mu_g}{\int_{\partial M}\sigma_{n+1-p}\,{\rm d}\mu_g}.
\end{equation}
Injecting (\ref{eq:constantsumnuintomegaI2}) again  into (\ref{eq:sumdeltaomegaI2bis}), we deduce that
\begin{equation}\label{eq:pointwiseequalityPoincareineqpformsboundaryinEuclideanspace}
\frac{n|H_0|^2-\frac{p-1}{n(n-1)}{\rm Scal}^{\partial M}}{\sigma_{n+1-p}}=\frac{\int_{\partial M}H\,{\rm d}\mu_g}{\int_{\partial M}\sigma_{n+1-p}\,{\rm d}\mu_g}\cdot\sigma_{n+1-p}.
\end{equation}
Note that (\ref{eq:pointwiseequalityPoincareineqpformsboundaryinEuclideanspace}) holds in any codimension~$m$.

Next we look at the space of parallel forms on $M$. Let $\hat{\omega}:=\left(\hat{\omega}_I\right)_I\in\bigoplus_I\Omega^p(M)$. Fixing a~pointwise o.n.b.\ $({\rm e}_j)_{1\leq j\leq n}$ of $T\partial M$, the family
\[\left\{{\rm e}_J,\nu\wedge {\rm e}_K \mid 1\leq j_1<\dots< j_p\leq n,\;1\leq k_1<\dots<k_{p-1}\leq n \right\}\]
is a pointwise o.n.b.\ of $\Lambda^pT^*M$. Decomposing each $\hat{\omega}_I$ in that pointwise basis and the canonical basis $({\rm d}x_I)_I$ of $\Lambda^p(\mathbb{R}^{n+m})^*$ (where the $p$-tuples $I$ are ordered) respectively allows us to consider~$\hat{\omega}$ as a pointwise matrix with $\binom{n+m}{p}$ rows and $\binom{n+1}{p}$ columns.
Note that, because the pointwise linear map $\iota_0^*\colon\Lambda^p(\mathbb{R}^{n+m})^*\to\Lambda^pT^*\partial M$ is  surjective, the $\bigl(\omega_I=\iota_0^*({\rm d}x_I)\bigr)_I$ obviously span $\Lambda^pT^*\partial M$, which already shows that the $\binom{n}{p}$ first columns of the matrix $\hat{\omega}$, namely $\bigl(\hat{\omega}_I({\rm e}_J)\bigr)_{I,J}$, must be linearly independent since that matrix has $\binom{n}{p}$ linearly independent rows. Next we would like to show that the rank of the whole matrix $\hat{\omega}$ is maximal, i.e., equal to $\binom{n+1}{p}$.

This already allows for finding expressions for the inner products of columns of the matrix~$\hat{\omega}$.
Namely, fix ${\rm e}_J$ and $\nu\wedge {\rm e}_K$ as above, then using equation \eqref{eq:nuinterioromega}, we compute
\begin{gather*}
\sum_I\hat{\omega}_I({\rm e}_J)\hat{\omega}_I(\nu\wedge {\rm e}_K)=\sum_I\langle\hat{\omega}_I,{\rm e}_J\rangle\cdot\langle\nu\lrcorner\,\hat{\omega}_I,{\rm e}_K\rangle\\
\hphantom{\sum_I\hat{\omega}_I({\rm e}_J)\hat{\omega}_I(\nu\wedge {\rm e}_K)}=-\frac{1}{\sigma_{n+1-p}}\!\sum_I\bigl\langle{\rm d}x_I,{\rm e}_J\bigr\rangle\biggl\langle\sum_{a=1}^m\mathbb{I}_{\nu_a}^{[p-1]}\bigl((\nu_a\lrcorner\,{\rm d}x_I)^T\bigr)-(nH_0\lrcorner\,{\rm d}x_I)^T\!, {\rm e}_K\biggr\rangle\\
\hphantom{\sum_I\hat{\omega}_I({\rm e}_J)\hat{\omega}_I(\nu\wedge {\rm e}_K)}=-\frac{1}{\sigma_{n+1-p}}\sum_I \bigl\langle {\rm d}x_I,{\rm e}_J\bigr\rangle\biggl\langle {\rm d}x_I,\sum_{a=1}^m\nu_a\wedge\mathbb{I}_{\nu_a}^{[p-1]}{\rm e}_K-nH_0\wedge {\rm e}_K\biggr\rangle\\
\hphantom{\sum_I\hat{\omega}_I({\rm e}_J)\hat{\omega}_I(\nu\wedge {\rm e}_K)}=-\frac{1}{\sigma_{n+1-p}}\biggl\langle {\rm e}_J,\sum_{a=1}^m \nu_a\wedge\mathbb{I}_{\nu_a}^{[p-1]}{\rm e}_K-nH_0\wedge {\rm e}_K \biggr\rangle
=0
\end{gather*}
because of both $\nu_a,H_0\perp T\partial M$. This shows that every among the $\binom{n}{p-1}$ last columns of $\hat{\omega}$, corresponding to the matrix $\bigl(\hat{\omega}_I(\nu\wedge {\rm e}_K)\bigr)_{I,K}$, is pointwise orthogonal to any of the $\binom{n}{p}$ first ones which correspond to the full-ranked matrix $\bigl(\hat{\omega}_I({\rm e}_J)\bigr)_{I,J}$.

We now look at the rank of the matrix $\bigl(\hat{\omega}_I(\nu\wedge {\rm e}_K)\bigr)_{I,K}$. For any $(p-1)$-tuples $J$, $K$, we compute
\begin{gather*}
\sum_I\hat{\omega}_I(\nu\wedge {\rm e}_J)\hat{\omega}_I(\nu\wedge {\rm e}_K)=\sum_I\langle\nu\lrcorner\,\hat{\omega}_I,{\rm e}_J\rangle\langle\nu\lrcorner\,\hat{\omega}_I,{\rm e}_K\rangle\\
 =\frac{1}{\sigma_{n+1-p}^2}\sum_I \biggl\langle {\rm d}x_I,\sum_{a=1}^m \nu_a\wedge\mathbb{I}_{\nu_a}^{[p-1]}{\rm e}_J-nH_0\wedge {\rm e}_J\biggr\rangle\cdot \biggl\langle {\rm d}x_I,\sum_{a=1}^m \nu_a\wedge\mathbb{I}_{\nu_a}^{[p-1]}{\rm e}_K-nH_0\wedge {\rm e}_K\biggr\rangle\\
 =\frac{1}{\sigma_{n+1-p}^2}\biggl\langle\sum_{a=1}^m \nu_a\wedge\mathbb{I}_{\nu_a}^{[p-1]}{\rm e}_J-nH_0\wedge {\rm e}_J,
\sum_{a=1}^m \nu_a\wedge\mathbb{I}_{\nu_a}^{[p-1]}{\rm e}_K-nH_0\wedge {\rm e}_K\biggr\rangle\\
 =\frac{1}{\sigma_{n+1-p}^2}\biggl(\sum_{a=1}^m\bigl\langle \mathbb{I}_{\nu_a}^{[p-1]}{\rm e}_J,\mathbb{I}_{\nu_a}^{[p-1]}{\rm e}_K\bigr\rangle
-2\bigl\langle\mathbb{I}_{nH_0}^{[p-1]}({\rm e}_J),{\rm e}_K\bigr\rangle+n^2|H_0|^2\bigl\langle {\rm e}_J,{\rm e}_K\bigr\rangle\biggr).
\end{gather*}
Here we notice that, in the particular case, where all $\mathbb{I}_{\nu_a}$ are simultaneously diagonalizable, i.e., if $[\mathbb{I}_{\nu_a},\mathbb{I}_{\nu_b}]=0$ for all $a$, $b$, we can choose $({\rm e}_j)_{1\leq j\leq n}$ so as to simultaneously diagonalize all~$\mathbb{I}_{\nu_a}$.
Then it is easy to show that $\sum_I\hat{\omega}_I(\nu\wedge {\rm e}_J)\hat{\omega}_I(\nu\wedge {\rm e}_K)=0$ for all $J\neq K$.
However, the $\sum_I\hat{\omega}_I(\nu\wedge {\rm e}_J)^2$ cannot be shown to be positive.
This means that, even if $[\mathbb{I}_{\nu_a},\mathbb{I}_{\nu_b}]=0$ for all~$a$,~$b$, one column of the matrix $\bigl(\hat{\omega}_I(\nu\wedge {\rm e}_K)\bigr)_{I,K}$ may vanish, in which case $\hat{\omega}$ will not be of full rank.\looseness=-1

From now on we furthermore assume that $\iota_0(\partial M)\subset V$, where $V$ is an $(n+1)$-dimensional affine subspace of $\mathbb{R}^{n+m}$.
Choosing a pointwise o.n.b.\ $(\nu_a)_{1\leq a\leq m}$ of $T^\perp\partial M\subset\mathbb{R}^{n+m}$ such that $\nu_1\in V$ and $\nu_a\perp V$ for all $a\geq2$, we obviously have $\mathbb{I}_{\nu_a}=0$ for all $a\geq2$ since $V\subset\mathbb{R}^{n+m}$ is totally geodesic. Therefore, choosing $({\rm e}_j)_{1\leq j\leq n}$ as an eigenbasis for the endomorphism $\mathbb{I}_{\nu_1}$ of~$T\partial M$, we have $\mathbb{I}_{\nu_1}{\rm e}_j=\kappa_j {\rm e}_j$ for all $1\leq j\leq n$ and the sum $\sum_I\hat{\omega}_I(\nu\wedge {\rm e}_J)\hat{\omega}_I(\nu\wedge {\rm e}_K)$ computed above simplifies to
\begin{gather*}
\sum_I\hat{\omega}_I(\nu\wedge {\rm e}_J)\hat{\omega}_I(\nu\wedge {\rm e}_K)\\
\qquad{}=\frac{1}{\sigma_{n+1-p}^2}\Big(\bigl\langle \mathbb{I}_{\nu_1}^{[p-1]}{\rm e}_J,\mathbb{I}_{\nu_1}^{[p-1]}{\rm e}_K\bigr\rangle
-2\langle\mathbb{I}_{nH_0}({\rm e}_J),{\rm e}_K\rangle+n^2|H_0|^2\langle {\rm e}_J,{\rm e}_K\rangle\Big)\\
\qquad{}=\frac{1}{\sigma_{n+1-p}^2}\biggl(\biggl(\sum_{j\in J}\kappa_j\biggr)^2-2\langle nH_0,\nu_1\rangle\biggl(\sum_{j\in J}\kappa_j\biggr)+\langle nH_0,\nu_1\rangle^2\biggr)\delta_{JK}\\
\qquad{}=\frac{1}{\sigma_{n+1-p}^2}\biggl(\sum_{j\in J}\kappa_j-\langle nH_0,\nu_1\rangle\biggr)^2\delta_{JK}
=\frac{1}{\sigma_{n+1-p}^2}\biggl(\sum_{j\notin J}\kappa_j\biggr)^2\delta_{JK},
\end{gather*}
where $\delta_{JK}=0$ if $J\neq K$ and $1$ if $J=K$. Now since $\iota_0\colon \partial M\to V$ is an isometric immersion of~an~$n$-dimensional manifold into an $(n+1)$-dimensional Euclidean space, there exists a point $x\in\partial M$ for which $\kappa_j(x)>0$ for all $1\leq j\leq n$ holds, see, e.g., \cite[p.~255]{GallotHulinLafontaine04}. At that point $x$, we can conclude that $\sum_I\hat{\omega}_I(\nu\wedge {\rm e}_J)\hat{\omega}_I(\nu\wedge {\rm e}_K)>0$ if $J=K$ and vanishes otherwise.
Therefore the columns of \smash{$\bigl(\hat{\omega}_I(\nu\wedge {\rm e}_K)\bigr)_{I,K}$ } form an orthogonal system of nonzero vectors at $x$, from which can be deduced that the whole matrix $\hat{\omega}$ has full rank $\binom{n+1}{p}$ at $x$. In turn, this implies that, at $x$, there are $\binom{n+1}{p}$ linearly independent rows in $\hat{\omega}$, that is $\binom{n+1}{p}$ linearly independent $\hat{\omega}_I$. Necessarily there must exist $\binom{n+1}{p}$ linearly independent parallel forms among the $\hat{\omega}_I$ on $M$, which is the maximal number allowed. As a first consequence, $\bigl(M^{n+1},g\bigr)$ must be flat (remember that $1\leq p\leq n$).
As a second consequence, at each point of $\partial M$, the family $(\nu\lrcorner\,\hat{\omega}_I)_I$ must span $\Lambda^{p-1}T^*\partial M$, so that $\iota\colon\partial M\to M$ must be totally umbilical by the identity $S^{[p-1]}(\nu\lrcorner\,\hat{\omega}_I)=(nH-\sigma_{n+1-p})\nu\lrcorner\,\hat{\omega}_I$ for all $I$ and the assumption $p\geq2$, see Remark \ref{r:LambdapTMtrivial}.
Since $M$ is flat and $\iota$ is totally umbilical in the Einstein manifold $M$, the mean curvature $H$ must be constant --  and positive because of $(n+1-p)H\geq\sigma_{n+1-p}>0$. Up to rescaling $g$ on $M$ as well as the Euclidean metric on $\mathbb{R}^{n+m}$, it may be assumed that $H=1$ along $\partial M$. By \cite[Theorem 13]{RaulotSavo11}, because $\bigl(M^{n+1},g\bigr)$ is flat, $\iota$ is totally umbilical and with constant mean curvature $1$, the manifold $\bigl(M^{n+1},g\bigr)$ must be isometric to the $(n+1)$-dimensional flat disk.
Moreover, identity (\ref{eq:pointwiseequalityPoincareineqpformsboundaryinEuclideanspace}) implies that $|H_0|=1=H$ along~$\partial M$. But then  by the proof of the equality case in \cite[Theorem 1.2]{MiaoWang14}, there exists  an isometric immersion $M\to V$ extending $\iota$. Since $\partial M\cong \mathbb{S}^n$ has constant sectional curvature $1$, the immersion $\iota_0$ must be an embedding by standard results due to Hadamard and Cohn--Vossen, see, e.g., \cite{doCarmoWarner70}. Again  by the proof of \cite[Theorem~1.2]{MiaoWang14},  the above isometric immersion $M\to V$ extending $\iota$  is  an~embedding. This shows that $\bigl(M^{n+1},g\bigr)$ is isometric to the $(n+1)$-dimensional flat disk standardly embedded in $V$. This concludes the proof of Theorem \ref{t:PoincareineqpformsboundaryinEuclideanspace}.
\end{proof}

\begin{Remarks}\quad \begin{enumerate}\itemsep=0pt
\item For $p=1$, inequality (\ref{eq:PoincareineqpformsboundaryinEuclideanspace}) reads
\begin{equation}\label{eq:ineqMiaoWang}
\int_{\partial M}\frac{|H_0|^2}{H}\,{\rm d}\mu_g\geq\int_{\partial M}H\,{\rm d}\mu_g
\end{equation}
assuming $\mathrm{Ric}=W^{[1]}\geq0$ on $M$ as well as $H=\frac{\sigma_n}{n}>0$ on $\partial M$. This is precisely the i\-ne\-qua\-li\-ty established by Miao and Wang in \cite[Theorem 1.2]{MiaoWang14} when $H>0$.
Actually, inequality~(\ref{eq:ineqMiaoWang})  implies (\ref{eq:PoincareineqpformsboundaryinEuclideanspace}) for all $1\leq p\leq n$ if not only $\sigma_{n+1-p}>0$ and $W^{[p]}\geq~0$ are assumed but also $\mathrm{Ric}\geq0$. The Gau\ss{} formula for curvature implies ${{\rm Scal}^{\partial M}=n^2|H_0|^2-|\mathbb{I}|^2}$ along $\partial M$. The  Cauchy--Schwarz i\-ne\-qua\-li\-ty yields $|\mathbb{I}|^2\geq n|H_0|^2$, so that ${\rm Scal}^{\partial M}\leq {n(n-1)|H_0|^2}$.
Fix now $p\in\{1,\dots,n\}$ and assume $\sigma_{n+1-p}>0$ along $\partial M$.
Because of $\frac{\sigma_{n+1-p}}{n+1-p}\leq H$, we can deduce that $H>0$ along $\partial M$ and that
\[\frac{n|H_0|^2-\frac{p-1}{n(n-1)}{\rm Scal}^{\partial M}}{\sigma_{n+1-p}}\geq\frac{n|H_0|^2-(p-1)|H_0|^2}{\sigma_{n+1-p}}\geq\frac{|H_0|^2}{H}\]
holds along $\partial M$.
Therefore, i\-ne\-qua\-li\-ty (\ref{eq:PoincareineqpformsboundaryinEuclideanspace}) can be deduced from i\-ne\-qua\-li\-ty (\ref{eq:ineqMiaoWang}) if $\mathrm{Ric}\geq0$ is also assumed. Mind however that our assumption $W^{[p]}\geq0$ differs from $\mathrm{Ric}\geq0$ for $2\leq p\leq n-1$, so that (\ref{eq:PoincareineqpformsboundaryinEuclideanspace}) cannot be deduced from \cite[Theorem~1.2]{MiaoWang14} in general.
\item According to \cite[Theorem 9]{RaulotSavo11}, given a compact Riemannian manifold $\bigl(M^{n+1},g\bigr)$ such that $W^{[p]}\geq 0$ for some $1\leq p\leq \frac{n+1}{2}$ and with boundary $\partial M$ isometric to the unit round sphere, then $M$ must be isometric to the Euclidean unit ball as soon as $\sigma_p\geq p$.
Actually this holds true for an \emph{arbitrary} $p\in\{1,\dots,n\}$. Namely if $\partial M$ is isometric to the round sphere~$\mathbb{S}^n$, then by taking $\iota_0\colon\mathbb{S}^n\to \mathbb{R}^{n+1}$ the standard embedding, we get, under the condition $\sigma_{n+1-p}>0$ together with the identities $H_0=1$ and ${\rm Scal}^{\partial M}=n(n-1)$
\[\int_{\partial M}\frac{n-p+1}{\sigma_{n-p+1}}\,{\rm d}\mu_g\geq \int_{\partial M}H\,{\rm d}\mu_g.\]
Therefore, if $\sigma_{n+1-p}\geq n+1-p$, then $H\geq1$ and the last i\-ne\-qua\-li\-ty is an equality, therefore~$M$ must be isometric to flat $\mathbb{D}^{n+1}$ by Theorem \ref{t:PoincareineqpformsboundaryinEuclideanspace}.
\end{enumerate}
\end{Remarks}

In the following, we consider the case where $W^{[p]}\geq p(n-p+1)\kappa$ for some nonvanishing real number $\kappa$.

\begin{Theorem}\label{thm:ineomegabiggerk}
Let $\bigl(M^{n+1},g\bigr)$ be any compact oriented $(n+1)$-dimensional Riemannian manifold with smooth nonempty boundary $\partial M$.
Fix $p\in\{1,\dots,n\}$ and assume that $W^{[p]}\geq p(n+1-p)\kappa$ on $M$ for some number $\kappa\neq0$ as well as $\sigma_{n+1-p}>0$ along $\partial M$.
Then, for any  $(p-1)$-form $\alpha$ on $\partial M$, we have
\begin{equation*}
\int_{\partial M}\frac{\bigl|\delta^{\partial M}{\rm d}^{\partial M}\alpha-\frac{p(n+1-p)\kappa}{2}\alpha\bigr|^2}{\sigma_{n+1-p}}\,{\rm d}\mu_g\geq\int_{\partial M}\bigl\langle S^{[p]}{\rm d}^{\partial M}\alpha,{\rm d}^{\partial M}\alpha\bigr\rangle\,{\rm d}\mu_g
\end{equation*}
with equality if and only if $\alpha=0$.
\end{Theorem}

\begin{proof} As in Theorem \ref{t:Poincareineqpformsgeneral}, we take the exact form $\omega={\rm d}^{\partial M}\alpha$ and consider the extension~$\hat\alpha$ on~$M$ such that $\delta {\rm d}\hat\alpha=0$ on $M$ with $\iota^*\hat\alpha=\alpha$ along $\partial M$.
The form $\hat\omega={\rm d}\hat\alpha$ satisfies ${\rm d}\hat\omega=0$, $\delta\hat\omega=0$ on $M$ and $\iota^*\hat\omega=\omega$.
We now apply identity (\ref{eq:Reillyformula}) to $\hat{\omega}$  to get after using that $\bigl\langle W^{[p]}\hat\omega,\hat\omega\bigr\rangle\geq p(n+1-p)\kappa|\hat\omega|^2$, inequality \eqref{eq:sn1p} and the fact $|\nabla\hat\omega|^2\geq 0$:
\[0\geq p(n+1-p)\kappa\int_M|\hat\omega|^2\,{\rm d}\mu_g+\int_{\partial M}\bigl(2\bigl\langle\nu\lrcorner\,\hat{\omega},\delta^{\partial M}\omega\bigr\rangle+\bigl\langle S^{[p]}\omega,\omega\bigr\rangle+\sigma_{n+1-p}|\nu\lrcorner\hat\omega|^2\bigr)\,{\rm d}\mu_g.
\]
Now, by Stokes formula, we have that
\[
\int_M|\hat\omega|^2\,{\rm d}\mu_g=\int_M\langle\hat\alpha,\underbrace{\delta {\rm d}\hat\alpha}_{0}\rangle\,{\rm d}\mu_g-\int_{\partial M}\langle\iota^*\hat\alpha,\nu\lrcorner {\rm d}\hat\alpha\rangle\,{\rm d}\mu_g=-\int_{\partial M}\langle\alpha,\nu\lrcorner\hat\omega\rangle\,{\rm d}\mu_g.
\]
Therefore, implementing this last equality into the previous i\-ne\-qua\-li\-ty yields
\[
0\geq \int_{\partial M}\left(2\left\langle\nu\lrcorner\,\hat{\omega},\delta^{\partial M}\omega-\frac{p(n+1-p)\kappa}{2}\alpha\right\rangle+\bigl\langle S^{[p]}\omega,\omega\bigr\rangle+\sigma_{n+1-p}|\nu\lrcorner\hat\omega|^2\right) {\rm d}\mu_g.
\]
By assumption $\sigma_{n+1-p}>0$ along $\partial M$, therefore the cross term
\[
t(\omega):=2\biggl\langle\nu\lrcorner\,\hat{\omega},\delta^{\partial M}\omega-\frac{p(n+1-p)\kappa}{2}\alpha\biggr\rangle
\] can be written as
\begin{align*} 
t(\omega)={}& \biggl|\sigma_{n+1-p}^{1/2}\nu\lrcorner\,\hat{\omega}+\frac{1}{\sigma_{n+1-p}^{1/2}}\bigl(\delta^{\partial M}\omega-\frac{p(n+1-p)\kappa}{2}\alpha\bigr)\biggr|^2\nonumber\\
&{}-\sigma_{n+1-p}|\nu\lrcorner\,\hat{\omega}|^2\nonumber-\frac{\bigl|\delta^{\partial M}\omega-\frac{p(n+1-p)\kappa}{2}\alpha\bigr|^2}{\sigma_{n+1-p}}\nonumber\\
 \geq{} & -\sigma_{n+1-p}|\nu\lrcorner\,\hat{\omega}|^2-\frac{\bigl|\delta^{\partial M}\omega-\frac{p(n+1-p)\kappa}{2}\alpha\bigr|^2}{\sigma_{n+1-p}}.
\end{align*}
Thus, we deduce after replacing $\omega$ by ${\rm d}^{\partial M}\alpha$ that
\[0\geq \int_{\partial M}\biggl(-\frac{\bigl|\delta^{\partial M}{\rm d}^{\partial M}\alpha-\frac{p(n+1-p)\kappa}{2}\alpha\bigr|^2}{\sigma_{n+1-p}}+\bigl\langle S^{[p]}{\rm d}^{\partial M}\alpha,{\rm d}^{\partial M}\alpha\bigr\rangle\biggr)\,{\rm d}\mu_g,
\]
which is the required i\-ne\-qua\-li\-ty. Notice that, if equality in the last i\-ne\-qua\-li\-ty occurs, then $\hat\omega$ is parallel, which implies $0=\langle W^{[p]}\hat\omega,\hat\omega\rangle= p(n+1-p)\kappa|\hat\omega|^2$ and thus $\hat\omega=0$. In turn this yields $\omega=0$ and, because of $p(n+1-p)\kappa\neq0$, also $\alpha=0$. This ends the proof.
\end{proof}

In the following, we will consider a compact Riemannian manifold $\bigl(M^{n+1},g\bigr)$ and assume furthermore that its boundary is immersed into the round sphere $\mathbb{S}^{n+m}(\kappa)$ of curvature $\kappa$.
\begin{Theorem}\label{t:MwithdelMinSn+m}
Let $\bigl(M^{n+1},g\bigr)$ be any compact oriented $(n+1)$-dimensional Riemannian manifold with smooth nonempty boundary $\partial M$. Assume that $W^{[p]}\geq p(n+1-p)\kappa$ on $M$ for some number $\kappa>0$ as well as $\sigma_{n+1-p}>0$ along $\partial M$ for a given $p\in\{1,\dots,n\}$.
Assume also that there exists an isometric immersion $\iota_0\colon\partial M\to\mathbb{S}^{n+m}(\kappa)$, with mean curvature vector $H_0$.
Then
\begin{equation*}
\int_{\partial M}\frac{n|H_0|^2-\frac{p-1}{n(n-1)}{\rm Scal}^{\partial M}+\frac{\bigl(np-(p-1)(p-4)\bigr)\kappa}{4}}{\sigma_{n+1-p}}\,{\rm d}\mu_g>\int_{\partial M}H\,{\rm d}\mu_g,
\end{equation*}
where, as before, ${\rm Scal}^{\partial M}$ denotes the scalar curvature of $\partial M$.
\end{Theorem}

\begin{proof} We will follow the same idea as in Theorem \ref{t:PoincareineqpformsboundaryinEuclideanspace}.
Let $\iota_1\colon\partial M\to \mathbb{R}^{n+m+1}$ be the isometric immersion with mean curvature vector $H_1$, which is the composition of the standard embedding $\mathbb{S}^{n+m}(\kappa)\hookrightarrow \mathbb{R}^{n+m+1}$ with $\iota_0$. Let us denote by ${\rm d}x_I:={\rm d}x_{i_1}\wedge\dots\wedge {\rm d}x_{i_p}\in\Lambda^p(\mathbb{R}^{n+m+1})^*$ and by $\omega_I:=\iota_1^* {\rm d}x_I\in\Omega^p(\partial M)$.
Let $\alpha_I$ the $(p-1)$-form on $\partial M$ given by $\alpha_I=x_{i_1}{}_{|_{\partial M}} {\rm d}x_{i_2}^T\wedge\dots\wedge {\rm d}x_{i_p}^T$, where ${\rm d}x_i^T=\iota_1^* {\rm d}x_i$.
Clearly, we have that $\omega_I={\rm d}^{\partial M}\alpha_I$. Now Theorem \ref{thm:ineomegabiggerk} applied to $\alpha_I$ gives that\looseness=-1
\begin{equation}\label{eq:inegualiteomegai}
\int_{\partial M}\frac{\bigl|\delta^{\partial M}\omega_I-\frac{p(n+1-p)\kappa}{2}\alpha_I\bigr|^2}{\sigma_{n+1-p}}\,{\rm d}\mu_g\geq\int_{\partial M}\bigl\langle S^{[p]}\omega_I,\omega_I\bigr\rangle\,{\rm d}\mu_g.
\end{equation}
Next, we want to sum \eqref{eq:inegualiteomegai} over $I$.
First the sum of the r.h.s over $I$ is equal to $p!\binom{n}{p}pH$ by~\eqref{eq:sumISpomegaI}. Now, for the l.h.s., we compute the sum
\[s:=\sum_I\bigl|\delta^{\partial M}\omega_I-\frac{p(n+1-p)\kappa}{2}\alpha_I\bigr|^2\]
as follows:
\begin{gather}
s=\sum_I\biggl(\bigl|\delta^{\partial M}\omega_I\bigr|^2+\frac{p^2(n+1-p)^2\kappa^2}{4}|\alpha_I|^2-p(n+1-p)\kappa\langle \delta^{\partial M}\omega_I,\alpha_I\rangle\biggr)\nonumber\\
\hphantom{}\stackrel{\eqref{eq:sumdeltaomegaI2}}{=}p!\binom{n}{p}p\biggl(n|H_1|^2-\frac{p-1}{n(n-1)}{\rm Scal}^{\partial M}\biggr)+\frac{p^2(n+1-p)^2\kappa}{4}\sum_{i_2,\dots,i_p}\bigl|{\rm d}x_{i_2}^T\wedge\dots\wedge {\rm d}x_{i_p}^T\bigr|^2\nonumber\\
\hphantom{s=}-p(n+1-p)\kappa\sum_I\langle \delta^{\partial M}\omega_I,\alpha_I\rangle. \label{eq:deltabord}
\end{gather}
Here $H_1$ is the mean curvature of $\iota_1$ which is related to the one of $\iota_0$ by $|H_1|^2=|H_0|^2+\kappa$.
Moreover,
\[
\sum_{i_1}^{n+m+1}x_{i_1}^2=\frac{1}{\kappa}
\]
since $\iota_0(\partial M)\subset\mathbb{S}^{n+m}(\kappa)$, sphere of radius $\kappa^{-1/2}$. In order to compute the last two sums in the above equality, we take $A={\rm Id}$ in \eqref{eq:traceap}, and thus, $A^{[p-1]}=(p-1) {\rm Id}$ and ${\rm tr}(A)=n$, to get that
\[\sum_{i_2,\dots,i_p}\bigl|{\rm d}x_{i_2}^T\wedge\dots\wedge {\rm d}x_{i_p}^T\bigr|^2=(p-1)!\binom{n}{p-1}.\]
On the other hand, denoting by $\mathbb{I}$ the second fundamental form of the immersion $\iota_1\colon\partial M\to \mathbb{R}^{n+m+1}$, we get from equation \eqref{eq:deltaomegaISavo} that
\begin{gather*}
\sum_I\langle \delta^{\partial M}\omega_I,\alpha_I\rangle = \sum_{i_1,\dots,i_p}x_{i_1}\biggl\langle\sum_{a=1}^{m+1}\mathbb{I}_{\nu_a}^{[p-1]}\bigl((\nu_a\lrcorner\,{\rm d}x_I)^T\bigr)-(nH_1\lrcorner\,{\rm d}x_I)^T,{\rm d}x_{i_2}^T\wedge\dots\wedge {\rm d}x_{i_p}^T\biggr\rangle\\
\qquad{} = -\kappa^{-1/2}\sum_{i_2,\dots,i_p}\biggl\langle\sum_{a=1}^{m+1}\mathbb{I}_{\nu_a}^{[p-1]}\bigl((\nu_a\lrcorner\,(\nu_1\wedge {\rm d}x_{i_2}\wedge\dots\wedge {\rm d}x_{i_p}))^T\bigr),{\rm d}x_{i_2}^T\wedge\dots\wedge {\rm d}x_{i_p}^T\biggr\rangle\\
\qquad\quad{} +\kappa^{-1/2}\bigl\langle\bigl(nH_1\lrcorner\,(\nu_1\wedge {\rm d}x_{i_2}\wedge\dots\wedge {\rm d}x_{i_p})\bigr)^T,{\rm d}x_{i_2}^T\wedge\dots\wedge {\rm d}x_{i_p}^T\bigr\rangle.
\end{gather*}
In the second equality, we use the fact that, in the local orthonormal basis $\{\nu_a\}_{a=1,\dots,m+1}$ of $T^\perp\partial M$ for the immersion $\iota_1\colon\partial M\to \mathbb{R}^{n+m+1}$, it may be assumed that $\nu_1=-\kappa^{1/2}\sum_{i_1}x_{i_1}dx_{i_1}$ which is the inner unit normal vector field for the standard immersion $\mathbb{S}^{n+m}(\kappa)\to\mathbb{R}^{n+m+1}$.
Hence, we proceed
\begin{align*}
\sum_I\bigl\langle \delta^{\partial M}\omega_I,\alpha_I\bigr\rangle ={}&-\kappa^{-1/2}\sum_{i_2,\dots,i_p}\bigl\langle\mathbb{I}_{\nu_1}^{[p-1]}({\rm d}x_{i_2}^T\wedge\dots\wedge {\rm d}x_{i_p}^T),{\rm d}x_{i_2}^T\wedge\dots\wedge {\rm d}x_{i_p}^T\bigr\rangle\\
&{} +\kappa^{-1/2}\langle nH_1,\nu_1\rangle\sum_{i_2,\dots,i_p} \bigl|{\rm d}x_{i_2}^T\wedge\dots\wedge {\rm d}x_{i_p}^T\bigr|^2\\
 \stackrel{\eqref{eq:traceap}}{=}{} &-\kappa^{-1/2}(p-1)!\binom{n}{p-1}\frac{p-1}{n}{\rm tr}(\mathbb{I}_{\nu_1})
+\kappa^{-1/2}\langle nH_1,\nu_1\rangle (p-1)!\binom{n}{p-1}\\
 ={} &(p-1)!\binom{n}{p-1}\frac{n-p+1}{n}\kappa^{-1/2}\langle nH_1,\nu_1\rangle.
\end{align*}
Now, the second fundamental form $\mathbb{I}$ of the immersion $\iota_1$ is the sum of the one of $\iota_0$ and the one of the isometric immersion $\mathbb{S}^{n+m}(\kappa)\to\mathbb{R}^{n+m+1}$.
Therefore, $\langle \mathbb{I}(X,Y),\nu_1\rangle=\kappa^{1/2} g(X,Y)$ for any $X,Y\in T\partial M$.
Thus, by tracing over an orthonormal frame of $T\partial M$, we get $\langle nH_1,\nu_1\rangle=n\kappa^{1/2}$. Then
\[\sum_I\bigl\langle \delta^{\partial M}\omega_I,\alpha_I\bigr\rangle=(p-1)!\binom{n}{p-1}(n-p+1).\]
Inserting this last computation into \eqref{eq:deltabord}, we finally deduce that
\begin{align*}
s={}&
p!\binom{n}{p}p\left(n|H_0|^2+n\kappa-\frac{p-1}{n(n-1)}{\rm Scal}^{\partial M}\right)\nonumber\\&{}+\frac{p^2(n-p+1)^2\kappa}{4}(p-1)!\binom{n}{p-1}-p!(n-p+1)^2\kappa\binom{n}{p-1}\\
={}& p!\binom{n}{p}p\left(n|H_0|^2+n\kappa-\frac{p-1}{n(n-1)}{\rm Scal}^{\partial M}\right)+(p-1)!\binom{n}{p-1}\frac{p (p-4 )(n-p+1)^2\kappa}{4}\\
={}& p!\binom{n}{p}p\left(n|H_0|^2+n\kappa-\frac{p-1}{n(n-1)}{\rm Scal}^{\partial M}+\frac{(p-4)(n-p+1)\kappa}{4}\right)\\
={}& p!\binom{n}{p}p\left(n|H_0|^2-\frac{p-1}{n(n-1)}{\rm Scal}^{\partial M}+\frac{(np-(p-1)(p-4))\kappa}{4}\right).
\end{align*}
Integrating the last identity and applying inequality \eqref{eq:inegualiteomegai} after simplifying by $p!\binom{n}{p}p$, we obtain the desired i\-ne\-qua\-li\-ty.
If equality holds, then for every multi-index $I$, necessarily $\alpha_I=0$ must hold by Theorem \ref{thm:ineomegabiggerk}.
But, pointwise, the $\alpha_I$'s span $\Lambda^{p-1}T^*\partial M$, therefore we obtain a~contradiction. This shows that the i\-ne\-qua\-li\-ty we obtained is strict and concludes the proof of~Theorem~\ref{t:MwithdelMinSn+m}.
\end{proof}

\section{A Ros-type i\-ne\-qua\-li\-ty for differential forms}\label{s:Rosdiffforms}

In this section, we generalize the Ros i\-ne\-qua\-li\-ty stated in \cite[Theorem 1]{Ros87} to the set-up of differential forms.
\begin{Theorem}\label{t:Rosgeneral}
Let $\bigl(M^{n+1},g\bigr)$ be any compact oriented $(n+1)$-dimensional Riemannian manifold with smooth nonempty boundary $\partial M$.
Fix $p\in\{2,\dots,n\}$  and assume that the $(p-1)^{\textrm{th}}$ relative de Rham cohomology group is reduced to zero, $W^{[p]}\geq0$ on $M$ as well as $\sigma_{n+1-p}>0$ along $\partial M$.
Assume also the existence of a nonzero parallel $(p-1)$-form $\omega_0$ on $M$, whose constant length may be assumed to be $1$.
Then
\begin{equation}\label{eq:ineqRosgeneral}
(n+2-p)\mathrm{Vol}(M,g)\leq(n+1-p)\int_{\partial M}\frac{|\iota^*\omega_0|^2}{\sigma_{n+1-p}}\,{\rm d}\mu_g.
\end{equation}
If equality in \eqref{eq:ineqRosgeneral} is realized, then under the assumptions that ${\rm Ric}\geq 0$ $($when $p\geq 3)$ and the mean curvature $H$ is constant, then the manifold $M$ is the $(n+1)$-dimensional flat disk $\mathbb{D}^{n+1}$.
\end{Theorem}

\begin{proof} Since we assume that $H_{\rm rel}^{p-1}(M)=0$, \cite[Theorem~3.2.5]{Schwarz95} implies the existence of a~$p$-form $\omega$ on $M$ with ${\rm d}\omega=0$, $\delta\omega=\omega_0$ on $M$ as well as $\iota^*\omega=0$ along $\partial M$. We apply (\ref{eq:Reillyformula}) to $\omega$.
 First, because ${\rm d}\omega=0$, it follows from \cite[Lemme 6.8]{GallotMeyer75} that
\begin{gather}\label{GalMey}
|\nabla\omega|^2\geq\frac{|\delta\omega|^2}{n-p+2},
\end{gather}
and since $|\delta\omega|^2=|\omega_0|^2=1$ and $W^{[p]}\geq0$ on $M$, we have
\begin{gather*}
\int_M\bigl(|{\rm d}\omega|^2+|\delta\omega|^2-|\nabla\omega|^2-\langle W^{[p]}\omega,\omega\rangle\bigr)\,{\rm d}\mu_g\\
\qquad{} \leq \int_M\Big(1-\frac{1}{n+2-p}\Big)\,{\rm d}\mu_g = \frac{n+1-p}{n+2-p}\mathrm{Vol}(M,g).
\end{gather*}
Moreover $\delta^{\partial M}(\iota^*\omega)=0$, $S^{[p]}\iota^*\omega=0$ and
\[\bigl\langle S^{[n+1-p]}\iota^*(\star\omega),\iota^*(\star\omega)\bigr\rangle\geq\sigma_{n+1-p}|\nu\lrcorner\,\omega|^2.\]
Therefore, (\ref{eq:Reillyformula}) yields
\[\frac{n+1-p}{n+2-p}\mathrm{Vol}(M,g)\geq\int_{\partial M}\sigma_{n+1-p}|\nu\lrcorner\,\omega|^2\,{\rm d}\mu_g.\]
Now partial integration together with Cauchy--Schwarz i\-ne\-qua\-li\-ty give
\begin{align}
\mathrm{Vol}(M,g)&= \int_M\langle\delta\omega,\omega_0\rangle\,{\rm d}\mu_g
=\int_M\langle\omega,\underbrace{{\rm d}\omega_0}_{0}\rangle\,{\rm d}\mu_g+\int_{\partial M}\langle\nu\lrcorner\,\omega,\iota^*\omega_0\rangle\,{\rm d}\mu_g\nonumber\\
&= \int_{\partial M}\biggl\langle\sigma_{n+1-p}^{1/2}\nu\lrcorner\,\omega,\frac{\iota^*\omega_0}{\sigma_{n+1-p}^{1/2}}\biggr\rangle\,{\rm d}\mu_g\nonumber\\
&\leq \biggl(\int_{\partial M}\sigma_{n+1-p}|\nu\lrcorner\,\omega|^2\,{\rm d}\mu_g\biggr)^{1/2}\cdot
\biggl(\int_{\partial M}\frac{|\iota^*\omega_0|^2}{\sigma_{n+1-p}}\,{\rm d}\mu_g\biggr)^{1/2},\label{eq:cauchyschwarz}
\end{align}
so that
\[\int_{\partial M}\sigma_{n+1-p}|\nu\lrcorner\,\omega|^2\,{\rm d}\mu_g\geq\frac{\mathrm{Vol}(M,g)^2}{\int_{\partial M}\frac{|\iota^*\omega_0|^2}{\sigma_{n+1-p}}\,{\rm d}\mu_g}.\]
Injecting that i\-ne\-qua\-li\-ty in the last one involving the volume of $M$, we obtain (\ref{eq:ineqRosgeneral}). Assume now equality in \eqref{eq:ineqRosgeneral} is attained, then equality holds in (\ref{GalMey}) so that
\[
\nabla_X\omega=-\frac{1}{n-p+2}X\wedge\delta\omega
\]
for all $X\in TM$. Also, we have that $S^{[n+1-p]}\iota^*(\star\omega)=\sigma_{n+1-p}\iota^*(\star\omega)$, which from the identity $*_{\partial M}S^{[p-1]}+S^{[n-p+1]}*_{\partial M}=nH*_{\partial M}$ valid on $(p-1)$-forms on $\partial M$ \cite[p. 624]{RaulotSavo11}, is equivalent to saying that $S^{[p-1]}(\nu\lrcorner\omega)=(nH-\sigma_{n+1-p})\nu\lrcorner\omega$. The Cauchy--Schwarz i\-ne\-qua\-li\-ty in \eqref{eq:cauchyschwarz} is an~equality and, thus, we get that $\iota^*\omega_0=c \sigma_{n+1-p}\nu\lrcorner\omega$ for some constant $c$. The constant $c$ can be determined by just replacing the last identity into \eqref{eq:ineqRosgeneral} and the second equality in \eqref{eq:cauchyschwarz} to deduce that $c=\frac{n+2-p}{n+1-p}$.
Now, from \cite[Lemma 18]{RaulotSavo11}, we get
\begin{equation}\label{eq:dbordnu}
{\rm d}^{\partial M}(\nu\lrcorner\omega)=-\nu\lrcorner {\rm d}\omega+\iota^*(\nabla_\nu\omega)-S^{[p]}\iota^*\omega=0.
\end{equation}
Here, we use that ${\rm d}\omega=0$, $S^{[p]}\iota^*\omega=0$ and that the tangential part of $\nabla_\nu\omega=-\frac{1}{n-p+2}\nu\wedge\delta\omega$ is zero as well.
On the other hand, we have that $\delta^{\partial M}(\nu\lrcorner\omega)=-\nu\lrcorner\delta\omega=-\nu\lrcorner\omega_0$. Hence, using~ \eqref{eq:dbordnu}, we get that
\begin{align*}
\Delta^{\partial M}(\nu\lrcorner\omega)&= -{\rm d}^{\partial M}(\nu\lrcorner\omega_0)=\nu\lrcorner {\rm d}\omega_0-\iota^*(\nabla_\nu\omega_0)+S^{[p-1]}\iota^*\omega_0\\
&= c\sigma_{n+1-p}(nH-\sigma_{n+1-p})\nu\lrcorner\omega.
\end{align*}
Taking now the $L^2$-scalar product of the last identity with $\nu\lrcorner\omega$ yields the following
\begin{align*}
c\int_{\partial M}\sigma_{n+1-p}(nH-\sigma_{n+1-p})|\nu\lrcorner\omega|^2&= \int_{\partial M}|\delta^{\partial M}(\nu\lrcorner\omega)|^2\,{\rm d}\mu_g
=\int_{\partial M}|\nu\lrcorner\omega_0|^2\,{\rm d}\mu_g\\
&=  {\rm Vol}(\partial M,g)-\int_{\partial M}|\iota^*\omega_0|^2\,{\rm d}\mu_g\\
&= {\rm Vol}(\partial M,g)-c^2\int_{\partial M}\sigma_{n+1-p}^2|\nu\lrcorner\omega|^2\,{\rm d}\mu_g.
\end{align*}
Hence, we deduce that
\[{\rm Vol}(\partial M,g)=c\int_{\partial M}\sigma_{n+1-p}\bigl(nH+(c-1)\sigma_{n+1-p}\bigr)|\nu\lrcorner\omega|^2\,{\rm d}\mu_g.\]
Now, the second equality in \eqref{eq:cauchyschwarz} gives that
\[{\rm Vol}(M,g)=c\int_{\partial M}\sigma_{n+1-p}|\nu\lrcorner\omega|^2\,{\rm d}\mu_g.\]
Hence, after replacing $c$ by $\frac{n+2-p}{n+1-p}$, we deduce that
\begin{align*}
\frac{{\rm Vol}(\partial M,g)}{{\rm Vol}(M,g)}&= \frac{\int_{\partial M}\sigma_{n+1-p}\bigl(nH+\frac{\sigma_{n+1-p}}{n+1-p}\bigr)|\nu\lrcorner\omega|^2\,{\rm d}\mu_g}{\int_{\partial M}\sigma_{n+1-p}|\nu\lrcorner\omega|^2\,{\rm d}\mu_g}\\
&\leq  \frac{\int_{\partial M}\sigma_{n+1-p}(nH+H)|\nu\lrcorner\omega|^2\,{\rm d}\mu_g}{\int_{\partial M}\sigma_{n+1-p}|\nu\lrcorner\omega|^2\,{\rm d}\mu_g}
= (n+1) H.
\end{align*}
The last equality uses the assumption that the mean curvature $H$ is constant. Therefore, we are in the equality case of Ros i\-ne\-qua\-li\-ty \cite[Theorem 1]{Ros87}. This allows to conclude the proof.
\end{proof}

Note that, when $p=1$, the assumptions of Theorem \ref{t:Rosgeneral} reduce to $\mathrm{Ric}\geq0$ on $M$ as well as $H>0$ along $\partial M$, in which case (\ref{eq:ineqRosgeneral}) reads
\[\int_{\partial M}\frac{1}{H}\,{\rm d}\mu_g\geq(n+1)\mathrm{Vol}(M,g)\]
which is exactly the Ros i\-ne\-qua\-li\-ty from \cite[Theorem 1]{Ros87}.

\begin{Remark}\rm Let us consider the particular case where $M^{n+1}$ is a domain in the Euclidean space $\mathbb{R}^{n+1}$. For any $p$, we take $i_1,\dots,i_{p-1}\in\{1,\dots,n+1\}$ such that $i_1<i_2<\dots<i_{p-1}$ and denote by $\omega_0:={\rm d}x_{i_1}\wedge\dots\wedge {\rm d}x_{i_{p-1}}$ the parallel $(p-1)$-form in $\Lambda^{p-1}\bigl(\mathbb{R}^{n+1}\bigr)^*$ which is of norm $1$. By Theorem \ref{t:Rosgeneral}, we have
\begin{equation*}
(n+2-p)\mathrm{Vol}(M,g)\leq(n+1-p)\int_{\partial M}\frac{\bigl|{\rm d}x_{i_1}^T\wedge\dots\wedge {\rm d}x_{i_{p-1}}^T\bigr|^2}{\sigma_{n+1-p}}\,{\rm d}\mu_g.
\end{equation*}
Summing over $i_1<i_2<\dots<i_{p-1}$, we get that
\begin{equation*}
\binom{n+1}{p-1}(n+2-p)\mathrm{Vol}(M,g)\leq(n+1-p)\int_{\partial M}\frac{\binom{n}{p-1}}{\sigma_{n+1-p}}\,{\rm d}\mu_g.
\end{equation*}
Here, we use that, by equation \eqref{eq:traceap} for $A={\rm Id}$,
\begin{gather*}
\sum_{i_1<i_2<\dots<i_{p-1}}\bigl|{\rm d}x_{i_1}^T\wedge\dots\wedge {\rm d}x_{i_{p-1}}^T\bigr|^2 = \frac{1}{(p-1)!}\sum_{i_1,i_2,\dots,i_{p-1}}\bigl|{\rm d}x_{i_1}^T\wedge\dots\wedge {\rm d}x_{i_{p-1}}^T\bigr|^2=\binom{n}{p-1}.
\end{gather*}
Hence, we deduce that
\begin{equation*}
(n+1)\mathrm{Vol}(M,g)\leq(n+1-p)\int_{\partial M}\frac{1}{\sigma_{n+1-p}}\,{\rm d}\mu_g
\end{equation*}
which actually can be deduced from Ros i\-ne\-qua\-li\-ty using $\frac{\sigma_{n+1-p}}{n+1-p}\leq H$.
\end{Remark}

\subsection*{Acknowledgements} We are very grateful to the \emph{Mathematisches For\-schungs\-in\-sti\-tut Oberwolfach} (MFO) and the \emph{Centre International de Rencontres Ma\-th\'e\-ma\-ti\-ques} (CIRM, Luminy) where most of the work was carried out.
The second-named author also thanks the Alfried Krupp Wissenschaftskolleg for its support.
{\it Last but not the least} we are grateful to the referees for their constructive comments.

\pdfbookmark[1]{References}{ref}
\LastPageEnding

\end{document}